\documentclass{commat}

\DeclareMathOperator{\ord}{ord}
\DeclareMathOperator{\supp}{supp}

\title{A family of non-Volterra quadratic operators corresponding to permutations}

\author{U. U. Jamilov}

\affiliation{
	Romanovskiy Institute of Mathematics, Uzbekistan Academy of Sciences. 9 University str., Tashkent, Uzbekistan 100174.  Akfa University. 264, National Park Street, Barkamol MFY, Yangiabad village, Qibray district, Tashkent, Uzbekistan 111221.  National University of Uzbekistan, 4,
University str., 100174, Tashkent, Uzbekistan.
	\email{uygun.jamilov@mathinst.uz, jamilovu@yandex.ru}
	}

\abstract{
In the present paper we consider a family of non-Volterra quadratic stochastic operators depending on a parameter $\alpha$ and study their trajectory behaviors.
We find all fixed points for a non-Volterra quadratic stochastic operator on a finite-dimensional simplex.  We construct some Lyapunov functions.
A complete description of the set of limit points is given, and we show that
such operators have the ergodic property.	
}



\VOLUME{31}
\YEAR{2023}
\NUMBER{1}
\firstpage{31}
\DOI{https://doi.org/10.46298/cm.10135}

\begin{paper}

\section{Introduction}

The quadratic stochastic operators frequently arise in many models of mathematical genetics, namely, in the theory of heredity
(see~\cite{ASZT}, \cite{Ber}, \cite{BlJaSc}, \cite{GGJ}, \cite{Gsb}, \cite{GMR}, \cite{GS}, \cite{Jlob}, \cite{Jjph}, \cite{Jdnc}, \cite{JL}, \cite{JLM}, \cite{JSW}, \cite{K}, \cite{Khukr}, \cite{L}, \cite{RZhmn}, \cite{RZhukr}, \cite{U}, \cite{Zus}, \cite{ZhRsb}).
Consider a~biological population and suppose that each individual in this population belongs precisely to one of the species (genotype)
 \(1,\dots, m\). The scale of species is such that the species of the parents \(i\) and \(j\),
unambiguously, determine the probability of every species \(k\) for the first generation of direct
descendants. Denote this probability, called the heredity coefficient, by \(p_{ij,k} = P(k|(i,j))\).
It is then obvious that \(p_{ij,k}\geq 0\) for all \(i,j,k\) and that
\[
\sum^m_{k = 1}p_{ij,k} = 1, \quad i,j,k = 1,\dots,m.
\]

The state of the population can be described by the tuple
\((x_1,x_2,\dots,x_m)\) of species probabilities, that is, \(x_k = P(k)\) is
the fraction of the species \(k\) in the total population. In the
case of panmixia (random interbreeding) the parent pairs \(i\) and
\(j\) arise for a~fixed state \(\mathbf{x} = (x_1,x_2,\dots,x_m)\) with
probability \(x_ix_j = P(i,j) = P(i)P(j)\). Hence, the total probability of the species
\(k\) in the first generation of direct descendants is defined by
\[
x'_k = \sum^m_{i,j = 1}P(k|(i,j)) P(i)P(j) = \sum^m_{i,j = 1}p_{ij,k}x_ix_j, \quad k = 1,\dots,m.
\]

The association \(\mathbf{x} \mapsto \mathbf{x}'\) defines an evolutionary quadratic operator.
 Thus evolution of a~population can be
studied as a~dynamical system of a~quadratic stochastic operator \cite{L}.
See~\cite{GMR} and~\cite{MG} for a~review of QSOs. Recently in~\cite{Khukr}, \cite{RHqtds} a~quasi-strictly non-Volterra QSO is studied.
We refer the reader to~\cite{JRjoca20} for a~review on convex combinations of quadratic stochastic operators.
The main goal of the present paper is to study a~family of operators which contains a~convex combination of two non-Volterra QSOs.

The paper is organised as follows. In Section~\ref{S:pre} we recall definitions and well known results from the theory of Volterra and non-Volterra QSOs.
In Section~\ref{S:mr} we consider a~class of non-Volterra QSOs and study trajectory behaviors of such operators.
We show that each QSO from this class has the two fixed points. Moreover, we prove that such operator
is~ergodic.

\section{Preliminaries}\label{S:pre}

Let
\[
    S^{m-1}
    = \left\{ \mathbf{x} = (x_1, x_2, \dotsc, x_m) \in \mathbb{R}^m \colon \textup{for any } i, \ x_i > 0 \textup{ and } \sum_{i = 1}^m x_i = 1 \right\}
\]
be the \((m-1)\)-dimensional simplex.
A map \(V\) of \(S^{m-1}\) into itself is called a~\textit{quadratic stochastic operator} (QSO) if
\begin{equation}\label{gqso}
(V \mathbf{x})_k = \sum_{i,j = 1}^m p_{ij,k}x_ix_j
\end{equation}
 for any \(\mathbf{x} \in S^{m-1}\) and for all \(k = 1,\dots,m\), where
\begin{equation}\label{condoper}
p_{ij,k} \ge 0, \quad
p_{ij,k} = p_{ji,k} \ \textup{ for all } i, j, k \quad \textup{and} \quad
\sum_{k = 1}^m p_{ij,k} = 1.
\end{equation}

Assume \(\{\mathbf{x}^{(n)}\in S^{m-1}: n = 0,1,2,\dots \}\) is the trajectory (orbit) of the initial point \(\mathbf{x}\in S^{m-1}\),
where \(\mathbf{x}^{(n+1)} = V(\mathbf{x}^{(n)})\) for all \(n = 0,1,2,\dots\), with \(\mathbf{x}^{(0)} = \mathbf{x}\).

One of the main problems in mathematical biology is to study the asymptotic behavior of the
trajectories. This problem deeply studied for the Volterra QSOs (see~\cite{Gsb}, \cite{Gmn}).
\begin{definition}
A quadratic stochastic operator is called a~\textit{Volterra operator} if
\begin{center}
\(p_{ij,k} = 0\) for any \(k\notin{{\{i,j\}}}, \ \ i, j, k = 1,\dots,m\).
\end{center}
\end{definition}

\begin{definition}
A point \(\mathbf{x}\in S^{m-1}\) is called a~\textit{periodic} point of \(V\) if there exists
an \(n\) so that \(V^n (\mathbf{x}) = \mathbf{x}\). The smallest positive integer \(n\) satisfying the above is called
the prime period or least period of the point \(\mathbf{x}\). A period-one point is called a~\textit{fixed} point of \(V\).
\end{definition}

Denote the set of all fixed points by \(\text{Fix}\, (V)\) and the set of all periodic points of
(not necessarily the smallest) period \(n\) by \(\text{Per}_n \, (V)\). Evidently that the set of all
iterates of a~periodic point form a~periodic trajectory (orbit).

Let \(D_\mathbf{x}V(\mathbf{x^*}) = (\partial V_i/\partial x_j)(\mathbf{x^*})\) be a~Jacobian of \(V\) at the point \(\mathbf{x^*}\).

\begin{definition}[\cite{Dev}]
A fixed point \(\mathbf{x}^*\) is called \textit{hyperbolic} if its Jacobian \(D_\mathbf{x}V(\mathbf{x}^*)\) has no eigenvalues on the unit circle.
\end{definition}

\begin{definition}[\cite{Dev}]
A hyperbolic fixed point \(\mathbf{x^*}\) is called:
\begin{itemize}
\item[i)] \textit{attracting}, if all the eigenvalues of the Jacobian \(D_\mathbf{x}V(\mathbf{x^*})\) are less than 1 in absolute value;
\item[ii)] \textit{repelling}, if all the eigenvalues of the Jacobian \(D_\mathbf{x}V(\mathbf{x^*})\) are greater than 1 in absolute value;
\item[iii)] a~\textit{saddle}, otherwise.
\end{itemize}
\end{definition}

\begin{definition}
A QSO \(V\) is called \textit{regular} if for any initial point \(\mathbf{x} \in S^{m-1}\), the limit
\(\lim\limits_{n\to \infty}V(\mathbf{x}^{(n)}) \) exists.
\end{definition}

Note that the limit point is a~fixed point of a~QSO. Thus, the fixed points
of a~QSO describe limit or long run behavior of the trajectories for any initial point.
The limit behavior of trajectories and fixed points play an important role in many
applied problems (see~\cite{BlJaSc}, \cite{GGJ}, \cite{Gsb}, \cite{GMR}, \cite{Jlob}, \cite{Jjph}, \cite{Juzmj}, \cite{JL}, \cite{JLM}, \cite{K}, \cite{L}, \cite{RZhmn}, \cite{Zus}, \cite{ZhRsb}).
The biological treatment of the regularity of a~QSO is rather clear: in the long run the
distribution of species in the next generation coincides with the distribution of species in the
previous one, i.e., it is stable.

For nonlinear dynamical systems~\eqref{gqso} Ulam~\cite{U} suggested an analogue of a~mea\-sure-\-theoretic ergodicity, the following ergodic hypothesis:

\begin{definition}
A QSO \(V\) is said to be \textit{ergodic} if the limit
\[
\lim_{n\to \infty} \frac{1}{n} \sum_{k = 0}^{n-1} V^k (\mathbf{x})
\]
exists for any \(\mathbf{x}\in S^{m-1}\).
\end{definition}

On the basis of numerical calculations Ulam, in~\cite{U}, conjectured that
the ergodic theorem holds for any QSO.
 In~\cite{Zus} Zakharevich proved that this conjecture is false in general.
 Later, in~\cite{GZ}, a~sufficient condition of non-ergodicity for QSOs defined on \(S^2\) was established.
 In~\cite{GGJ} have shown the correlation between non-ergodicity of Volterra QSOs
 and rock-paper-scissors games. In~\cite{JSW} the random
dynamics of Volterra QSOs is studied.

 The biological
treatment of non-ergodicity of a~QSO is the following: in the long run the behavior of the distributions of species is unpredictable.
Note that a~regular QSO is ergodic, but in general from ergodicity does not follow regularity.

Let \(\omega_V\big(\mathbf{x}^{(0)}\big)\) be the set of limit points of the trajectory
\[
\left\{ V^n \left( \mathbf{x}^{(0)} \right) \in S^{m-1} \colon n = 0,1,2, \dotsc \right\}.
\]

\begin{definition}
A continuous function \(\varphi \colon S^{m-1}\rightarrow R\) is called a~Lyapunov function for a~QSO \(V\) if
\(\varphi(V (\mathbf{x}))\geq \varphi(\mathbf{x})\) for all \(\mathbf{x}\) (or \(\varphi(V (\mathbf{x}))\leq \varphi(\mathbf{x})\) for
all \(\mathbf{x}\)).
\end{definition}

Note that a~Lyapunov function is very helpful to describe an upper estimate of \(\omega_V(\mathbf{x}^0)\).

\begin{definition}
A permutation \(\pi\) of \(E_n = \{1,\dots, n\}\) is a~\(k\)--\textit{cycle} if there exists a~positive integer \(k\)
and an integer \(i\in E_n\) such that
\begin{itemize}
\item[(1)] \(k\) is the smallest positive integer such that \(\pi^k (i) = i\), and
\item[(2)] \(\pi\) fixes each \(j\in E_n \setminus \{i, \pi(i), \dots , \pi^{k-1}(i)\}\).
\end{itemize}
The \(k\)-cycle \(\pi\) is usually denoted \(\big(i, \pi(i), \dots , \pi^{k-1}(i)\big)\).
\end{definition}

The set \(\supp (\pi) = \{i\in E_n: \pi(i)\neq i \}\) denotes the support of \(\pi\) and we let \(\supp (k)\) denote the support of the \(k\)-cycle, that is, the set
\[
    \supp (k) = \{i, \pi(i), \dots , \pi^{k-1}(i)\}.
\]

Any permutation can be represented in the form of a~product of cycles without common elements (i.e. disjoint cycles) and this representation is unique to within the order of the factors.

Let \(\pi = \tau_1\tau_2\dots\tau_q\) be a~permutation of the set \(E_{m-1} = \{1,\dots,m-1\}\), where \(\tau_1,\dots,\tau_q\) are disjoint cycles and we denote by \(\ord (\tau_i)\) the order of a~cycle \(\tau_i\). Evidently that
\[
\supp(\tau_1)\cup\dots\cup \supp(\tau_q) = \supp (\pi) \ \ \text{and} \ \ \supp(\tau_i)\cap \supp(\tau_j) = \emptyset, \ \ \text{for any} \ \ i\neq j.
\]

The following notations will be used in the below. Let
\[
\partial S^{m-1}
= \left\{ \mathbf{x} \in S^{m-1} \colon x_i = 0 \textup{ for at least one } i \in \{ 1, 2, \dotsc, m\} \right\}
\]
denote the boundary of \(S^{m-1}\) and let \(\text{int}\, S^{m-1} = \left\{ \mathbf{x} \in S^{m-1} \colon x_1 x_2 \dotsm x_m > 0 \right\}\) be the interior of \(S^{m-1}\).

\section{Main results}\label{S:mr}

Consider a~non-Volterra QSO defined on a~finite-dimensional simplex which has the form
 \begin{equation}\label{aop}
V_{\pi}: \left\{
\begin{array}
{ll}
 x'_k = 2 x_m x_{\pi(k)}, \ \ k = 1,\dots, m-1 \\
 x'_m = x_m ^2 +\Big(\sum\limits_{i = 1}^{m-1} x_i\Big)^2 \\
 \end{array}
\right.
\end{equation}
where \(\pi\) is a~permutation on the set \(E_{m-1}\).

It is worth mentioning that if \(\pi = (21)(3)\) then the QSO~\eqref{aop} coincides up to the rearrangement of the coordinates
with the quasi-strictly non-Volterra QSO which is studied in~\cite{Khukr}.

Let \(s = \text{LCM} \, \big(\ord (\tau_1),\dots, \ord (\tau_q)\big)\).
\begin{theorem}[\cite{JLKhsaa20}] \label{pertr} For the operator \(V_{\pi}\) the following statements are true:
\begin{itemize}
\item[i)] if \(\mathbf{x}^{(0)}\in \Gamma = \{\mathbf{x}\in S^{m-1}: x_m = 0\}\cup \{\mathbf{e}_m\}\) then
\(\omega_{V_{\pi}}\big(\mathbf{x}^{(0)}\big) = \{\mathbf{e}_m\}\);
\item[ii)] if \(\pi = Id\) then \(\omega_{V_{\pi}}\big(\mathbf{x}^{(0)}\big) = \{\widetilde{\mathbf{x}} \}\) for any \(\mathbf{x}^{(0)}\in S^{m-1}\setminus \Gamma\);
\item[iii)] if \(\pi\neq Id\) then \(\omega_{V_{\pi}}\big(\mathbf{x}^{(0)}\big) = \{\mathbf{x}_\xi, \mathbf{x}_\xi^{1},\dots, \mathbf{x}_\xi^{s-1}\}\).
\end{itemize}
\end{theorem}

Let \(V_1 := V_{Id}\) and \(V_2 := V_{\pi}\). Consider the convex combination of the QSOs \(V_1,V_2\), that~is,
\[
V_\alpha = \alpha V_1 +(1-\alpha)V_2, \ \ \alpha\in [0,1]
\]
It is easy to see that the operator \(V_\alpha\) has the form
\begin{equation}\label{vop}
V_\alpha: \left\{
\begin{array}
{ll}
 x'_k = 2 x_m (\alpha x_k+ (1-\alpha) x_{\pi(k)}), \ \ k = 1,\dots, m-1 \\
 x'_m = x_m ^2 +\Big(\sum\limits_{i = 1}^{m-1} x_i\Big)^2 \\
 \end{array}
\right.
\end{equation}
where \(\pi\) is a~permutation on the set \(E_{m-1}\).

It is evident that if \(\pi = Id\) then for any \(\alpha\in [0,1]\) the operator \(V_\alpha\)
coincides with the QSO \(V_1\). The dynamics of the operator \(V_1\) is given in the
Theorem~\ref{pertr}. In the below we consider the cases \(\pi\neq Id\).

The QSO~\eqref{vop} can be written as follow
\begin{equation}\label{vop1}
V_\alpha: \left\{
\begin{array}
{lll}
 x'_k = 2 x_m (\alpha x_k+ (1-\alpha) x_{\pi(k)}), \ \, \ \ k\in \supp(\pi) \\
 x'_k = 2 x_m x_k, \ \ k\notin \supp(\pi)\\
 x'_m = x_m ^2 +\Big(\sum\limits_{i = 1}^{m-1} x_i\Big)^2 \\
 \end{array}
\right.
\end{equation}
where \(\pi\) is a~permutation on the set \(E_{m-1}\).

Consider the function
\(f(x) = 2x^2 -2x+1, \ \ x\in [0,1]\).
We define \(f^n\) to be the \(n\)-fold composition of \(f\) with itself.
One can easily verify the statements of the next proposition about dynamics of the function \(f(x)\).
\begin{proposition}\label{f} For the function \(f(x)\) the following statements are true:
\begin{itemize}
\item[i)] \(\text{Fix} \, (f) = \{1, 1/2\}\);
\item[ii)] \(x = 1\) is a~repelling and the fixed point \(1/2\) is an attracting;
\item[iii)] for any value \(n\geq2\) the function \(f(x)\) has no \(n-\) periodic points, different from fixed points;
\item[iv)] \(\lim\limits_{n\rightarrow\infty} f^n (x) = 1/2\) for any \(0<x<1\) and \(f(0) = f(1) = 1\).
\end{itemize}
\end{proposition}

Denote \(\supp (\mathbf{x}) = \{i: \, x_i>0\}\) and let \(|\supp (\mathbf{x})|\) be its cardinality.

 In the next Proposition we will describe the invariant sets, all fixed points and
we give some Lyapunov function.
\begin{proposition}\label{proper}
For the operator \(V_\alpha\) the following statements are true:
\begin{itemize}
\item[i)] If \(|\supp \, (\pi)|<m-1\) then
\(\Gamma_\beta = \{\mathbf{x}\in S^{m-1}: x_i = 0, \, \forall \, i\in\beta\}\) is an~invariant set for any \(\beta\subset E_{m-1}\setminus \supp \, (\pi)\). Also the sets
\[
M_{\mu,i} = \bigg\{\mathbf{x}\in S^{m-1}: \, \sum_{k\in \supp(\tau_i)} x_{k} = \mu, \ \ \, x_m = 1/2 \bigg\} \ \ \text{and}
\]
 \[
M_{\nu,i,j} = \bigg\{\mathbf{x}\in S^{m-1}: \, \sum_{k\in \supp(\tau_i)} x_{k} = \nu \, \sum_{k\in \supp(\tau_j)} x_k \,\bigg\}
\]
are invariant sets, where \(\mu\geq0, \nu >0\);
\item[ii)] \(\text{Fix}\, (V_\alpha) = X\cup \{\mathbf{e}_m\}\), where \(\mathbf{e}_m = (0,\dots,0,1)\) and
 \[
X = \{\mathbf{x}\in S^{m-1}: \, x_k = x_l, \,\ \forall \, k,l\in\supp(\tau_i), \ \ i = 1,\dots,q, \ \ x_m = 1/2\};
\]
\item[iii)] For any \(i\in \{1,\dots,q\}\) the function \(\varphi_i\big(\mathbf{x}\big) = \sum\limits_{k\in \supp(\tau_i)} x_{k}\)
 is a~Lyapunov function;
\item[iv)] For any \(k\notin\supp (\pi)\) the function \(\phi_k\big(\mathbf{x}\big) = x_k \) is a~Lyapunov function.
\end{itemize}
\end{proposition}

\begin{proof}
i) Let \(|\supp \, (\pi)|<m-1\) then for any \(k\notin \supp \, (\pi)\) from~\eqref{vop1} one easily has \(x'_k = 0\). Hence it follows that the set
\(\Gamma_\beta\) is a~invariant set.

Let \(\mathbf{x} \in M_{\mu,i}\) and \(\tau_i\) is a~cycle then from~\eqref{aop} we have
\begin{align*}
\sum\limits_{k\in \supp(\tau_i)} x'_{k}
&= \sum\limits_{k\in \supp(\tau_i)} \big(\alpha x_{k}+(1-\alpha)x_{\pi(k)}\big) \\
&= \alpha \sum\limits_{k\in \supp(\tau_i)} x_{k} + (1-\alpha) \sum\limits_{k\in \supp(\tau_i)} x_{\pi(k)} \\
&= \mu .
\end{align*}
Therefore \(V(M_{\mu,i})\subset M_{\mu,i}\).

Let \(\mathbf{x} \in M_{\nu,i,j}\) and \(\tau_i,\tau_j\) cycles then from~\eqref{aop} we have
\[
\frac{\sum\limits_{k\in \supp(\tau_i)} x'_{k}}{\sum\limits_{k\in \supp(\tau_j)} x'_{k}} =
\frac{2x_m\sum\limits_{k\in \supp(\tau_i)} x_{k}}{2x_m\sum\limits_{k\in \supp(\tau_j)} x_{k}} = \nu.
\]
Consequently \(V(M_{\nu,i,j})\subset M_{\nu,i,j}\).

ii) The equation \(V_\alpha(\mathbf{x}) = \mathbf{x}\) has the following form
\begin{equation}\label{eqaop}
\left\{
\begin{array}
{llllll}
 x_{k} = 2x_m (\alpha x_k+(1-\alpha) x_{\pi(k)}), \ \ 1\leq k\leq m-1,\\
[2mm]
 x_{m} = 2x^2_m-2x_m+1.
\end{array}
\right.
\end{equation}

Due to Proposition~\ref{f} the last equation of the system~\eqref{eqaop} has
the solutions \(x_m = 1\) and \(x_m = 1/2\).

Evidently that if \(x_m = 1\) then we get the vertex \(\mathbf{e}_m = (0,\dots,0,1)\).

For \(x_m = 1/2\) from the system of equations
\[
x_{k} = \alpha x_k+(1-\alpha) x_{\pi(k)}, \ \ 1\leq k\leq m-1, \ \ \text{and} \ \ x_1+\dots+x_{m-1} = \frac{1}{2}
\]
it follows that
\[
x_k = x_{k'} \ \ \text{for all} \ \ k,k'\in \supp (\tau_i), \ \ i = 1,\dots,q \ \ \text{and}
\]
\[
x_k = x_k \ \ \text{for all} \ \ k\in \text{Fix} (\pi).
\]

Using the last one has that a~point \(\mathbf{x} = (x_1,\dots,x_m)\in X\) is a~solution of the system~\eqref{eqaop}.

iii) Let \(\tau_i, i\in\{1,\dots, q\}\) be a~cycle. Then \(f(x)\geq 1/2\) for any \(0<x<1\) and we can assume that
\(x_m\geq1/2\). Then from~\eqref{aop} we have

\begin{align*}
\varphi_i\big(V_\alpha(\mathbf{x})\big) &= \sum\limits_{k\in \supp(\tau_i)} x'_{k} = \sum\limits_{k\in \supp(\tau_i)} 2x_m (\alpha x_k+(1-\alpha)x_{\pi(k)})\\
&=
2x_m \Big(\alpha\sum\limits_{k\in \supp(\tau_i)} x_k +(1-\alpha) \sum\limits_{k\in \supp(\tau_i)} x_{\pi(k)} \Big)\\
&= 2x_m \Big(\alpha \varphi_i(\mathbf{x}) +(1-\alpha)\varphi_i(\mathbf{x})\Big) = 2x_m\varphi_i(\mathbf{x})\\
&\geq \varphi_i(\mathbf{x}).
\end{align*}

Therefore the functions \(\varphi_i(\mathbf{x})\) are Lyapunov functions for any \(i\in\{1,\dots, q\}\).

iv) Using \(x_m\geq 1/2\) from~\eqref{vop1} for any \(k\notin\supp(\pi)\) one easily has that
\[
\phi_k(V_\alpha(\mathbf{x})) = 2x_mx_k\geq x_k = \phi_k(\mathbf{x}).
\qedhere
\]
\end{proof}

\begin{corollary}
If \(p = |E_{m-1}\setminus \supp \, (\pi)|\), then
\[
\phi(\mathbf{x}) = \gamma_1\phi_1(\mathbf{x}) + \dots + \gamma_p \phi_p(\mathbf{x}) \ \ \text{and} \ \
\varphi(\mathbf{x}) = \beta_1\varphi_1(\mathbf{x}) + \dots + \beta_q \varphi_q(\mathbf{x})
\]
are Lyapunov function for the QSO \(V_\alpha\) for any \(\gamma_1 \ge 0, \dotsc, \gamma_p \ge 0\) and \(\beta_1 \ge 0, \dotsc, \beta_q \ge 0\).
\end{corollary}

 In the next Theorem we give the description of the set of limit points of the trajectories.

\begin{theorem}\label{pertrcc} For the operator \(V_\alpha\) the following statements are true:
\begin{itemize}
\item[i)] if \(\mathbf{x}^{(0)}\in \Gamma = \{\mathbf{x}\in S^{m-1}: x_m = 0\}\cup \{\mathbf{e}_m\}\) then
\(\omega_{V_\alpha}\big(\mathbf{x}^{(0)}\big) = \{\mathbf{e}_m\}\);
\item[ii)] if \(\alpha\in(0,1), \ \ \pi\neq Id\) then \(\omega_{V_\alpha}\big(\mathbf{x}^{(0)}\big) = \{\mathbf{b} \}, \ \ \mathbf{b}\in X\) for any \(\mathbf{x}^{(0)}\in S^{m-1}\setminus (\Gamma\cup X)\);
\item[iii)] if \(\alpha = 0, \ \ \pi\neq Id\) then \(\omega_{V_\alpha}\big(\mathbf{x}^{(0)}\big) = \{\mathbf{x}_\xi, \mathbf{x}_\xi^{1},\dots, \mathbf{x}_\xi^{s-1}\}\);
\item[iv)] if \(\alpha = 1, \ \ \pi\neq Id\) then \(\omega_{V_\alpha}\big(\mathbf{x}^{(0)}\big) = \{\widetilde{\mathbf{x}} \}\) for any \(\mathbf{x}^{(0)}\in S^{m-1}\setminus \Gamma\).
\end{itemize}
\end{theorem}

\begin{proof}
i) Evidently that \(V_\alpha\big(\mathbf{x}^{(0)}\big) = \mathbf{e}_m\) for any \(\mathbf{x}^{(0)}\in \Gamma\).

ii) Let \(\alpha\in(0,1)\) and \(\mathbf{x}^{(0)}\in S^{m-1}\setminus (\Gamma\cup X)\).
 Then by assertion of Proposition~\ref{f} we obtain \(\lim\limits_{n\rightarrow\infty} x_m^{(n)} = 1/2\).
 Since \(f(x)\geq 1/2\) for any \(0<x<1\) and we can assume that \(x_m\geq1/2\).

 Let \(k\notin \supp \, (\pi)\).
Due to Proposition~\ref{proper} the function \(\phi_k\big(\mathbf{x}\big)\) is a~Lyapunov function for the QSO~\eqref{vop}.
Therefore, we have
\begin{equation}\label{112}
\phi_k(\mathbf{x}^{(n+1)})\geq \phi_k(\mathbf{x}^{(n)}), \ \ k\notin\supp \, (\pi), \ \ n = 0,1,\dots,
\end{equation}
that is there exists \(\lim\limits_{n\rightarrow\infty} x^{(n)}_k = \lim\limits_{n\rightarrow\infty}\phi_k(\mathbf{x}^{(n)}) = \xi_k\) for any \(k\notin\supp \, (\pi)\).

Denote \(\widetilde{X} = \Big\{\mathbf{x}\in S^{m-1}: b_k = \xi_k, \, \forall \, k\notin\supp \, (\pi), \, b_m = 1/2\Big\}\).

Let \(\tau_i, i\in\{1,\dots, q\}\) be a~cycle. Consider the function \(\psi_i(\mathbf{x}) = \min\limits_{k\in\supp(\tau_i)} x_k\).
Then from~\eqref{aop} we have
\begin{equation}\label{psibah}
\psi_i\big(V_\alpha(\mathbf{x})\big) = \min\limits_{k\in\supp(\tau_i)} 2x_m (\alpha x_k+(1-\alpha)x_{\pi(k)}) \geq \alpha \psi_i(\mathbf{x})+ (1-\alpha)\psi_i(\mathbf{x}) = \psi_i(\mathbf{x}).
\end{equation}
 Consequently, we have
\begin{equation}\label{12}
\psi_i(\mathbf{x}^{(n+1)})\geq \psi_i(\mathbf{x}^{(n)}), \ \ i = 1,\dots,q, \ \ n = 0,1,\dots
\end{equation}
Therefore the sequence \(\big\{\psi_i(\mathbf{x}^{(n)})\big\}\) is an increasing and bounded sequence.
Hence it follows existence the following limit \(\lim\limits_{n\rightarrow\infty}\psi_i(\mathbf{x}^{(n)}) = \xi_i\).

Let \(k\in\supp \, (\pi)\). It is easy to see that for any \(i\in\{1,\dots,q\}\)
\[
\psi_i\big(\mathbf{x}\big)\leq \psi_i\big(\mathbf{b}\big) \ \ \text{and} \ \
\psi_i\big(\mathbf{x}\big) = \psi_i\big(\mathbf{b}\big) \ \ \text{iff} \ \ \mathbf{x} = \mathbf{b}, \ \ \mathbf{b}\in \widetilde{X}\cap X.
\]
Indeed if \(\mathbf{x} = \mathbf{b}, \ \ \mathbf{b}\in \widetilde{X}\cap X\) then it is easily follows that \(\psi_i\big(\mathbf{x}\big) = \psi_i\big(\mathbf{b}\big)\) for any \(i\in\{1,\dots,q\}\). Let \(\psi_i\big(\mathbf{x}\big) = \psi_i\big(\mathbf{b}\big), \, \mathbf{x}\in \widetilde{X}\) for any \(i\in\{1,\dots,q\}\) and \(\mathbf{b}\in \widetilde{X}\cap X\) then for any \(i\in\{1,\dots,q\}\) we get \(\psi_i\big(\mathbf{x}\big) = x_{i_1}\leq x_{i_2}\leq \dots \leq x_{i_t}\), where \(t = \ord (\tau_i)\). If we assume that some of the inequalities \(x_{i_1}\leq x_{i_2}\leq \dots \leq x_{i_t}\) are strong inequalities in this case we have contradiction to \(\mathbf{x}\in S^{m-1}\). Therefore we have if \(\psi_i\big(\mathbf{x}\big) = \psi_i\big(\mathbf{b}\big)\) for any \(i\in\{1,\dots,q\}\) and \(\mathbf{b}\in \widetilde{X}\cap X\) then for any \(i\in\{1,\dots,q\}\) we obtain \(\psi_i\big(\mathbf{x}\big) = x_{i_1} = x_{i_2} = \dots = x_{i_t}\).

Next we prove that if \(\xi_i<\psi_i(\mathbf{b})\) for any \(i\in\{1,2,\dots,q\}\),
then \(\lim\limits_{n\rightarrow\infty} \mathbf{x}^{(n)} = \mathbf{b}\). Suppose the converse. Then there is a~sequence
\(\{\mathbf{x}^{(n_t)}\}_{t = 1,2,3,\dots}\) such that
\begin{equation}\label{limnt}
\lim\limits_{t\rightarrow\infty} \mathbf{x}^{(n_t)} = \mathbf{c}\neq\mathbf{b}.
\end{equation}
Using \(\text{min} f(x) = 1/2\) one has
\begin{align*}
1& = \frac{\psi_i(\mathbf{b})-\xi_i}{\psi_i(\mathbf{b})-\xi_i} = \lim\limits_{t\rightarrow\infty}\frac{\psi_i(\mathbf{b})-\psi_i(\mathbf{x}^{(n_t+1)})}{\psi_i(\mathbf{b})-\psi_i(\mathbf{x}^{(n_t)})}\\
&= 1+ \lim\limits_{t\rightarrow\infty} \frac{\psi_i(\mathbf{x}^{(n_t)})-2x_m^{(n_t)}\Big(\alpha x_{k}^{(n_t)}+(1-\alpha)x_{\pi(k)}^{(n_t)}\Big)}{\psi_i(\mathbf{b})-\psi_i(\mathbf{x}^{(n_t)})}\\
&\leq 1+ \lim\limits_{t\rightarrow\infty} \frac{\psi_i(\mathbf{x}^{(n_t)})- \Big(\alpha x_{k}^{(n_t)}+(1-\alpha)x_{\pi(k)}^{(n_t)}\Big)}{\psi_i(\mathbf{b})-\psi_i(\mathbf{x}^{(n_t)})}\\
&\leq 1+ \lim\limits_{t\rightarrow\infty} \frac{\psi_i(\mathbf{x}^{(n_t)})-1}{\psi_i(\mathbf{b})-\psi_i(\mathbf{x}^{(n_t)})}\\
&<1.
\end{align*}
This is a~contradiction. It follows that \(\xi_i = \psi_i(\mathbf{b})\) for any \(i\in\{1,2,\dots,q\}\).

Thus \(\lim\limits_{n\rightarrow\infty} \mathbf{x}^{(n)} = \mathbf{b}\) for any \(\alpha\in(0,1)\) and an initial \(\mathbf{x}^{(0)}\in S^{m-1}\setminus (\Gamma\cup X)\).

The proofs of parts iii) and iv) follows from the Theorem~\ref{pertr}.
\end{proof}

\begin{corollary}
The QSO \(V_\alpha\) is an ergodic transformation.
\end{corollary}

\subsection*{Acknowledgments}
The author thanks the referees for useful comments which helped to improve the presentation.

\EditInfo{April 13, 2020}{September 09, 2020}{Utkir Rozikov}

\end{paper}